\documentclass[11pt]{article}
\pdfoutput=1
\usepackage{amsthm}
\usepackage{latexsym, amssymb, amsmath}
\usepackage{graphicx}
\usepackage[dvipsnames]{xcolor}
\usepackage{fancyvrb}

\begin{document}
\numberwithin{equation}{section}

\newtheorem{THEOREM}{Theorem}
\newtheorem{PRO}{Proposition}
\newtheorem{XXXX}{\underline{Theorem}}
\newtheorem{CLAIM}{Claim}
\newtheorem{COR}{Corollary}
\newtheorem{LEMMA}{\underline{Lemma}}
\newtheorem{REM}{Remark}
\newtheorem{EX}{Example}
\newenvironment{PROOF}{{\bf Proof}.}{{\ \vrule height7pt width4pt depth1pt} \par \vspace{2ex} }
\newcommand{\Bibitem}[1]{\bibitem{#1} \ifnum\thelabelflag=1 
  \marginpar{\vspace{0.6\baselineskip}\hspace{-1.08\textwidth}\fbox{\rm#1}}
  \fi}
\newcounter{labelflag} \setcounter{labelflag}{0}
\newcommand{\labelon}{\setcounter{labelflag}{1}}
\newcommand{\Label}[1]{\label{#1} \ifnum\thelabelflag=1 
  \ifmmode  \makebox[0in][l]{\qquad\fbox{\rm#1}}
  \else\marginpar{\vspace{0.7\baselineskip}\hspace{-1.15\textwidth}\fbox{\rm#1}}
  \fi \fi}

\newcommand{\LEFTLINE}{\ifhmode\newline\else\noindent\fi}
\newcommand{\RIGHTLINE}[1]{\LEFTLINE\rightline{#1}}
\newcommand{\CENTERLINE}[1]{\LEFTLINE\centerline{#1}}
\def\BOX #1 #2 {\framebox[#1in]{\parbox{#1in}{\vspace{#2in}}}}
\parskip=8pt plus 2pt
\def\AUTHOR#1{\author{#1} \maketitle}
\def\Title#1{\begin{center}  \Large\bf #1 \end{center}  \vskip 1ex }
\def\Author#1{\vspace*{-2ex}\begin{center} #1 \end{center}  
 \vskip 2ex \par}
\renewcommand{\theequation}{\arabic{section}.\arabic{equation}}
\def\bdk#1{\makebox[0pt][l]{#1}\hspace*{0.03ex}\makebox[0pt][l]{#1}\hspace*{0.03ex}\makebox[0pt][l]{#1}\hspace*{0.03ex}\makebox[0pt][l]{#1}\mbox{#1} }
\def\psbx#1 #2 {\mbox{\psfig{file=#1,height=#2}}}

 
\newcommand{\ia}{\,\,\nearrow}
\newcommand{\da}{\,\,\searrow}
\newcommand{\ic}{\nearrow}
\newcommand{\dc}{\searrow}
\newcommand{\me}{\mbox{e}}
\newcommand{\vp}[2]{\vphantom{\vrule height#1pt depth#2pt}}
\renewcommand{\thefootnote}{\fnsymbol{footnote}}
\newcommand{\sr}{{\par\centerline{\vrule height0.1pt width2.5in depth0pt}\par}}
\newcommand{\fb}[1]{\fbox{ $\displaystyle #1 $ }}
\newcommand{\FG}[2]{\fbox{\begin{minipage}{5.4in}\vspace*{5pt}\includegraphics[width=#1mm]{#2.JPG}\vspace*{5pt}\end{minipage}}}
\newcommand{\fg}[2]{{\includegraphics[height=#1mm]{#2.JPG}}}
\newcommand{\F}{{\stepcounter{figure}\thefigure}}
 
 
\Title{Technical Details of the Proof of the Sine Inequality \\[1.2ex] {\normalsize  $\displaystyle  \sum_{k=1}^{n-1}\left( \frac{n}{k} - \frac{k}{n} \right) ^\beta  \sin(kx) \geq 0$ }} 


\begin{center}
MAN KAM KWONG
\end{center}

\begin{center}
\emph{Department of Applied Mathematics\\ The Hong Kong Polytechnic University,\\ Hunghom, Hong Kong}\\
\tt{mankwong@polyu.edu.hk}
\end{center}

\par\vspace*{\baselineskip}\par

\newcommand{\Cr}{\color{red}}
\newcommand{\mb}{\mathbf}
\newcommand{\tx}[1]{\mbox{ #1 }}

\parskip=6pt

\CENTERLINE{\underline{\em This is a preliminary draft, temporarily used as a place holder.}}
\CENTERLINE{\underline{A final version will be provided as soon as possible.}}

\par\vspace*{\baselineskip}\par
\begin{abstract}
In a recent study, H. Alzer and the author showed that
the sine polynomial
$$  \sum_{k=1}^{n-1} \left( \frac{n}{k} - \frac{k}{n} \right) ^\beta  \,\sin(kx) > 0   $$
is nonnegative for $ x\in[0,\pi ] $,
$ n\geq 2, \, \beta \geq  \beta _1 := \frac{\log(2)}{\log(16/5)} . $
This result, among others, will be presented in a forthcoming
article.
The proof relies on quite a number of technical Lemmas and inequalities.
We have decided to delegate all the tedious details of the proofs
of these Lemmas in a separate article, namely, the current one.
Some of the proofs require brute-force numerical computation,
performed with the help of the computer software MAPLE.

A few of the Lemmas included here are of independent interest.
\end{abstract}

\vspace{0.9cm}

\vspace{0.2cm}

\section{Introduction }
As stated in the Abstract, this article is meant to be a supplement to 
a hopefully forthcoming paper on a recent project done in collaboration
with H. Alzer. The results presented here are proved mainly by 
elementary techniques, sometimes with the assistance
of the computer for carrying out tedious computations.
There is no intention to seek formal publication in an official research
journal. For that reason, often more detailed arguments as well as
heuristic discussions may be included.

The result to be established is

\vspace{0.3cm}
\noindent
{\bf{Theorem AK.}} \emph{The inequality}
\begin{equation}
0\leq \sum_{k=1}^{n-1} \Bigl(\frac{n}{k}-\frac{k}{n}\Bigr)^{\beta}\sin(kx)
\quad{(\beta\in\mathbb{R})}  \Label{sp1}
\end{equation}
\emph{holds for all integers $n\geq 2$ and real numbers $x\in [0,\pi]$ if and only if} $\beta \geq \log(2)/\log(16/5)=0.59592...$.

\par\vspace*{1mm}\par
Throughout the paper, $ a_k,b_k,c_k,d_k $ denote NN (nonnegative) numbers, and $ m,n,k $ positive
integers.
For convenience, we adopt the notations:
$$  s(k) = \sin(kx), \qquad  c(k) = \cos(kx), \hspace*{5mm}  $$
\begin{eqnarray*}
         [a_1,a_2, ..., a_n] &=& \sum_{k=1}^{n} a_k \, s(k) ,  \\
     {}   [a_1,a_2, ..., a_n]^- &=& \sum_{k=1}^{n} (-1)^{k+1}a_k \, s(k) .
\end{eqnarray*}

The following identities are well-known:
\begin{equation}  \tau _k \,\, = \,\,\, [\underbrace{1, ... \, , 1\vp{0}{5}}_{k \tx{\footnotesize  terms}}]  \,\, =  \,\, \frac{c\big(\frac{1}2\big)-c\big({k}+\frac12\big)}{2s\big(\frac12\big)} .  \hspace*{9mm}  \Label{tau}  \end{equation}
\begin{equation}  \tau _k^-  \, =  \, [\, 1, ... \, , 1 \,]^-  \, =  \, \frac{s\big(\frac{1}2\big)-(-1)^{k}s\big({k}+\frac12\big)}{2c\big(\frac12\big)} \,.  \Label{tau2}  \end{equation}

\par\vspace*{\baselineskip}\par
A tool used frequently below is explained in detail in \cite{k1}.
Let $ f(x)=g_1(x)-g_2(x) $ be the difference of two
specific increasing (or decreasing) functions in a given bounded interval $ [\alpha ,\beta ] $ and they do 
not involve any further variable parameters.
The tool, called {\tt dif}, is a numerical algorithm, implemented 
as a MAPLE procedure, can be applied to 
rigorously prove $ f(x)\geq 0 $ (if that is true).
In the sequel, if an inequality is proved using this technique, we will attach the
computer output.

The author would like to thank Horst Alzer for an enjoyable 
long-time collaboration in the study of classical inequalities. The
current work is a part of one of our projects. Horst should be a 
coauthor, but he insisted otherwise. 

\section{Proof of Theorem AK}

The result for small $ n=2,\cdots ,6 $ was proved in \cite{ak}. Thus our starting
point is $ n\geq 7 $.

The proof, in its entirety, is quite long. In order not to
mask the main ideas, we postpone the proofs, some rather technical, of
quite a number of Lemmas and inequalities to Section 3. In addition, we divide the
proof into small chunks in order to highlight the different ideas involved.

{
\addtolength{\leftmargini}{-3mm}
\begin{enumerate}
\item Instead of studying the sine polynomial in (\ref{sp1}) directly, we consider 
two equivalent ones.
First, we divide it
by the coefficient of the first term to get
\begin{equation}  S_{n,\beta }(x) := [ a_{n,k} ] = \sum_{k=1}^{n-1} \left(  \frac{n^2-k^2}{(n^2-1)k} \right) ^\beta  s(k)   \Label{sn}  \end{equation}
with leading coefficient $ a_{n,1}=1 $.
For convenience, the dependence of $ a_{n,k} $ on $ \beta  $ is not explicit in
the notation. In addition, we
often suppress even $ n $ and{/}or $ x $ in the notation.
For example, we write $ S $ and sometimes $ S_n $
instead of $ S_{n,\beta }(x) $ and $ a_k $ instead of $ a_{n,k} $.
Note that $ a_k>0 $ for $ 1\leq k<n $, while $ a_n=0 $.

We also denote the sequence of second differences of the coefficients by 
\begin{equation}  \Box_k = a_{k-1} - 2 a_k + a_{k+1} \qquad  (k=2,\cdots ,n-1).  \end{equation}

Next, by reflection, we see that $ [a_{k}] $ is NN if and only if 
$ S^-=[a_{k}]^- $ is NN. (\underline{$\vphantom{y}$Caution}: In \cite{ak}, the superscript $ ^* $ is 
used instead of $ ^- $.)
The $ ^- $ operator can be extended to general sine polynomials.
Using $ S^- $ in placed of $ S $ simplifies some
intermediate inequalities. 

\item Our first step is to show that for $ n\geq 7 $, $ S^- $ is NN on $ [0.75,\pi ] $.

For all $ n $,
$ a_{2}\leq 2^{-\beta _1}<0.662 $. Using (\ref{cc1}) of Lemma~\ref{ck}, we have
\begin{eqnarray*}
   { [ a_{k} ]^-} \geq  \frac{ \left[ s\big(\frac{1}2\big)+s\big(\frac{3}2\big)\right] - 0.662\left[ 1+s\big(\frac{3}2\big)\right] }{2c\big(\frac{1}2\big)} =  \frac{ s\big(\frac{1}2\big)+0.338s\big(\frac{3}2\big)- 0.662}{2c\big(\frac{1}2\big)} \,.
\end{eqnarray*}
The numerator of the RHS is NN in $ [0.75,\pi ] $ (\cite[Lemma 7]{k1}).
so is $ [a_k]^- $.

\item A sequence $ \left\{ c_k\right\} _{k=1}^m $ is said to be convex (concave) if 
$$  c_{k-1} - 2 c_k + c_{k+1} \geq  (\,\leq \,) \,\, 0 \qquad  \tx{for}  k=2,\cdots ,m-1.  $$

Lemma~\ref{Fej} (ii) assets that when the coefficient sequence of a sine polynomial, appended with a 0,
is convex, the polynomial is NN.
For certain combinations of $ n $ and $ \beta  $ (for instance, when $ n=7 $ and 
$ \beta >0.8714) $,
the coefficient sequence $ \left\{ a_k\right\} _{k=1}^n $ (Note that the last member of the sequence is
$ a_n=0 $.) of (\ref{sn}) is convex, and hence $ S $
is NN, and Theorem~AK hold for those combinations.

However, in general, $ \left\{ a_k\right\}  $ is not convex. In such cases, we invoke Lemma~\ref{fc}.
to find a positive integer $ m<n $ such that the sub-sequence $ \left\{ a_1,\cdots ,a_m\right\}  $ is
convex while the sub-sequence $ \left\{ a_m,a_{m+1},\cdots ,a_n\right\}  $ has an odd number of terms
and is concave. After subtracting $ a_m $ from each of the terms of the 
first sub-sequence, we see that $ \left\{ a_1-a_m,\cdots ,a_{m-1}-a_m,0\right\}  $ is convex.
By Lemma~\ref{Fej} (ii), 
$$  S_1 := [\, \underbrace{a_1-a_m,\,a_2-a_m,\,\cdots ,\,a_{m-1}-a_m}_{m-1\tx{\footnotesize  terms}} \, ]  $$
is NN.
Note that $ S $ has the decomposition
\begin{eqnarray*}
   S &=& S_1 + T \\
     &=& S_1 + [\, \underbrace{a_{m}, \cdots ,a_{m}\vp{0}{5}}_{m\tx{\footnotesize  terms}}, a_{m+1}, ... , a_{n-1}\, ] .
\end{eqnarray*}

It is easy to verify that $ T $ has the alternative representation
\begin{equation}  T = \sum_{k=m}^{n-1} d_k \tau _k ,     \Label{Tdk}  \end{equation}
with an even number of summands on the RHS and
$$  d_k=a_k-a_{k+1} \ic \tx{in} k .  $$

\item By Lemma~\ref{fc} (ii), $ S_1 $ has at least five terms. We decompose it further as
follows
$$  S_1 = H + K  $$
where
\begin{equation}  H = [ \, a_1-4a_4+3a_5, \,\, a_2-3a_4+2a_5, \,\, a_3-2a_4+a_5 \, ] .  \Label{SH}  \end{equation}
The convexity of $ \left\{ a_1,\cdots ,a_m\right\}  $ implies that the coefficients of $ H $ and $ K $ are 
positive and satisfy the hypotheses of Lemma~\ref{Fej} (ii), respectively. Hence, both $ H $ and $ K $ are NN.
It follows that
\begin{equation}  S = H + K + T  \geq  H + T .  \Label{SHT}  \end{equation}
In view of this, an appropriate tactic to prove 
Theorem~AK is to find good lower bounds for $ H $ and $ T $.

\item For convenience, we have suppressed the explicit dependence of $ H $ and $ T $
on $ n $ and $ \beta  $.  The more precise notations ought to be
$ H(n,\beta ) $ and $ T(n,\beta ) $.

After the substitution $ y=\frac{1}{n^2-1} \in(0,\frac{1}{48} ) $,
the coefficients of $ H $ are given by $ h_i(y;\beta ) $, $ i=1,2,3 $
while the two second differences are $ h_i(y;\beta ) $, $ =4,5 $
as given before Lemma~\ref{hhh}.

For fixed $ n $ and $ \beta >\beta _1 $, Lemma~\ref{hhh} (ii) asserts that $ H(n,\beta )-H(n,\beta _1) $
satisfies Lemma~\ref{Fej} (ii) and so is NN, implying that $ H(n,\beta )\geq H(n,\beta _1) $.

Likewise, with $ \beta =\beta _1 $ fixed, and $ n>n_1 $, Lemma~\ref{hhh} (i) asserts that
$ H(n,\beta _1)\geq H(\nu _1,\beta _1) $.

It follows that for $ n\geq n_1 $ and $ \beta \geq \beta _1 $,
\begin{equation}  H(n,\beta )\geq H(n,\beta _1)\geq H(n_1,\beta _1).  \Label{Hnb1}  \end{equation}

\item In particular, for $ n>7 $ and $ \beta >\beta _1 $,
\begin{eqnarray*}
     && H^-(n,\beta ) \\
     &\geq &H^-(7,\beta _1) \\
     &=& \textstyle  \left( 1 - 4 \left( \frac{11}{64} \right) ^{\beta _1} + 3 \left( \frac{1}{10} \right) ^{\beta _1} \right) \sin(x) - \left(  \left( \frac{15}{32} \right) ^{\beta _1}- 3 \left( \frac{11}{64} \right) ^{\beta _1} + 2 \left( \frac{1}{10} \right) ^{\beta _1}  \right) \sin(2x) \\
     && \textstyle  \hspace*{3mm}  + \left( \left( \frac{5}{18} \right) ^{\beta _1} - 2 \left( \frac{11}{64} \right) ^{\beta _1} +  \left( \frac{1}{10} \right) ^{\beta _1} \right)  \sin(3x) .
\end{eqnarray*}
By Lemma~\ref{H7b}, this sine polynomial
is concave in $ [0,0.75] $, and so $ H^-(7,\beta _1)/x $ is decreasing, implying that
\begin{equation}  \frac{H^-}{x} \geq  \left . \frac{H^-(7,\beta _1)}{x} \, \right |_{x=0.75} = 0.2232352723...  \Label{H1}  \end{equation}

\item To get a lower bound for $ T^- $ we make use of (\ref{Tdk}) and Lemma~\ref{T3} (iii).
By definition,
$$  d_{n-1}=a_{n-1} = \left(  \frac{2n-1}{(n^2-1)(n-1)} \right) ^{\beta } \leq  \left(  \frac{2n-1}{(n^2-1)(n-1)} \right) ^{\beta _1} \leq  \, \frac{1.105}{n} .   $$
The last inequality holds because 
$ n \left(  \frac{2n-1}{(n^2-1)(n-1)} \right) ^{\beta _1} $
is a decreasing function of $ n $ (\cite[Lemma 13]{k1})
and the upper bound $ 1.105 $ is obtained by letting $ n=7 $.
Then Lemma~\ref{T3} (iii) gives
\begin{eqnarray}
         T^- &\geq & \frac{1.105}{n} \min\left\{ \tau _{n-1}^-,0\right\}  \nonumber  \\
              &\geq & \frac{1.105}{n} \left(  \frac{s\big(\frac12\big)-1}{2c\big(\frac12\big)} \right)  \nonumber  \\
	      &\geq & - \, \frac{1.105}{2n\cos(0.375)} \, .   \Label{T4}
\end{eqnarray}
\item From (\ref{SHT}), (\ref{H1}), and (\ref{T4}), we get
$$  S^- \geq   H^- +  T^- \geq  0.2232x - \frac{1.105}{2n\cos(0.375)} \, .  $$
The righthand side is NN when
\begin{eqnarray}
   x &\geq & \left(  \frac{1.105}{2(0.2232)\cos(0.375)} \right)  \frac{1}{n} \nonumber  \\[1.5ex]
     &>& \frac{2.660223693}{n} \,.   \Label{xn2}
\end{eqnarray}
This greatly improves \P2.

\item Now, (\ref{H1}) can be improved, 
by using the shorter interval $ [0,\frac{2.67}{7} ] $
instead of $ [0,0.75] $, to give
\begin{equation}  \frac{H^-}{x} \geq  \left . \frac{H^-(7,\beta _1)}{x} \, \right |_{x=2.67/7} > 0..2285 .  \Label{H2}  \end{equation}
This can be used to replace $ 0.2232 $ in the first line 
of (\ref{xn2}). 

Then, the factor $ \cos(0.375) $ in the denominator
of the same expression can be replaced by the larger number
$ \cos\big(\frac{2.67}{2(7)} \big) $. After that, we conclude that
$ S^- $ is NN when
\begin{eqnarray*}
  x &\geq & \left(  \frac{1.105}{2(0.2285)\cos\big(\frac{2.67}{2(7)}\big)} \right)  \frac{1}{n} \\
     &>& \frac{2.4602482}{n}  .
\end{eqnarray*}
This proves {\bf Proposition 1} of \cite{ak}.

Although the same arguments can be used iteratively to bootstrap the
assertion further, but the improvements gained this way are very slight.

\item When $ n $ is even, (\ref{tau2}) gives
$$  \tau _{n-1}^- = \frac{s\big(\frac{1}2\big)+s\big(n-\frac12\big)}{2c\big(\frac12\big)} = \frac{s\big(\frac{n}2\big)c\big(\frac{n-1}2\big)}{c\big(\frac12\big)} \geq  0, \qquad  x \in \textstyle  [0,\frac{\pi }{n} ] \,.   $$
By Lemma~\ref{T3} (iii), $ T^- $ is NN in  
$ [0,\frac{\pi }{n} ]\supset[0,\frac{2.5}{n} ] $. By (\ref{SHT}), 
$ S^- $ is thus also NN in $ [0,\frac{\pi }{n} ] $. Combining with \P9 yields Theorem~AK.
This is {\bf Proposition 2} of \cite{ak}.

It now remains to show that $ S^- $ is NN in $ [0,\frac{2.5}{n} ] $ for odd $ n\geq 7 $.

\item Let us determine $ T $ for $ n=7 $. The last second difference
$$  \Box_6 = \left( \frac{1}{10} \right) ^\beta  -2 \left( \frac{13}{288} \right) ^\beta   $$
is nonnegative for $ \beta \geq \beta _2:=\frac{\ln(2)}{\ln(288)-\ln(130)} \approx 0.8714162659 $. For such
values of $ \beta  $, $ T=0 $ and Theorem~AK holds. On the other hand,
$$  \Box_5 = \left( \frac{11}{64} \right) ^\beta  - 2 \left( \frac{1}{10} \right) ^\beta  + \left( \frac{13}{288} \right) ^\beta  > 0  $$
for all $ \beta >\beta _1 $. The construction presented in \P2 yields $ T $ with two
summands in the form of $ Tdk $.

\item For general $ n $, we can find an upper bound of the number of summands of $ T $ as
follows. We first compute $ \Box_k $ for $ S $ when $ \beta =\beta _1 $. Suppose
that the last $ N $ differences are negative, then $ T $ is constructed using the
last $ N $ or $ N-1 $ coefficients, whichever is even. By Remark \ref{cov}, for
$ \beta >\beta _1 $, $ T $ has fewer or equal number of summands than that for $ \beta _1 $.

Following this scheme, we find that for $ n\leq 13 $, $ T $ has no more
than 2 summands while for $ n\leq 43 $, there are no more than 10.

Recall that in all cases, $ T $ has an even number of summands in (\ref{Tdk})
and (\ref{d2}) of Lemma~\ref{T3} (iv) is applicable to yield a lower bound.

In (\ref{d2}),
$$  d_{2j}-d_{2j-1} = 2a_{2j}-a_{2j-1}-a_{2j+1}= - \Box_{2j}.  $$
By defining
$ \displaystyle \delta _k := - \Box_{n-k}, $
(\ref{d2}) becomes
\begin{equation}  \frac{2|T^-|}{x} \leq   (m+1)\delta _{n-m-1} + ... + (n-3)\delta _{3} + (n-1)\delta _{1} .  \Label{Tdm}  \end{equation}
We only need to worry about $ x\in[0,\frac{2.5}{n} ] $ in which $ T^- <0 $;
this explains the use of $ |T^-| $ on the LHS. Note that
only odd subscripts of $ \delta  $ are involved. 
If $ T $ has $ 2j $ summands, the
RHS of (\ref{Tdm}) has $ j $ terms involving $ \delta _1,\delta _3,\cdots ,\delta _{2j-1} $. The necessary
lower bounds for $ \delta _k $ are derived in Lemmas 12 and 13.

\item Let us first treat the simplest case when
  $T$ has 2 summands. This is true for $ n\leq 13 $.

Applying (\ref{gh1}) of Lemma 12 to (\ref{Tdm}) containing only one summand,
leads to
$$  \frac{|T^-|}{x} \leq  0.196 .  $$
With (\ref{H2}), this implies $ S^- \geq  H^- + T^- \geq 0 $ in $ [0,\frac{2.5}{n} ] $
and Theorem~AK is proved.

\item Next, suppose $T$ has 10 or less summands. This is true for $ n\leq 43 $.

For the rest of the proof we can assume $ n\geq 15 $ odd. By (\ref{Hnb1}),
we have $ H(n,\beta )\geq H(15,\beta _1) $.
In addition, we use the shorter interval $ [0,\frac{2.5}{15} ] $ to
improve (\ref{H2}) to
\begin{equation}  \frac{H^-}{x} \geq  \left . \frac{H^-(15,\beta _1)}{x} \, \right |_{x=2.5/15} > 0..248 \,.  \Label{H3}  \end{equation}

Using (\ref{Tdm}) with $ \delta _k $ up to $ \delta _{9} $ and
the estimates (\ref{gh9}) and (\ref{H3}), we obtain
$$  \frac{|T^-|}{x} \leq  0.196+0.0206+0.009+0.005171+0.003451 = 0.234222 < \frac{H^-}{x} \,,  $$
proving the Theorem.

\item Finally we consider the case when $T$ has more than 10 summands.

Using $ H(n,\beta )\geq H(45,\beta _1) $ and the interval $ \big[0,\frac{2.5}{45} \big] $,
we improve (\ref{H3}) further to
\begin{equation}  \frac{H^-}{x} \geq  \left . \frac{H^-(45,\beta _1)}{x} \, \right |_{x=2.5/45} > 0.250772629 \,.  \Label{H4}  \end{equation}

From (\ref{Tdm}), we obtain
\begin{equation}  \frac{2|T^-|}{x} \leq   (n-11) \big[\, \delta _{n-m-1} + ... + \delta _{11} \,\big] + (n-9) \delta _{9} + ... + (n-1) \delta _{1} .  \Label{Tdm2}  \end{equation}
The monotonicity of $ \Box_k $ implies 
$$  \delta _{n-m-1} + ... + \delta _{11} \leq  \delta _{n-m} + ... + \delta _{10} .  $$
Note that the subscripts of $ \delta  $ on the RHS are even; those 
on the LHS are odd. It follows from (\ref{Tdm2}) that
$$  \frac{2|T^-|}{x} \leq   \frac{(n-11)}{2} \, \big[\, \delta _{n-m-1} + ... + \delta _{10} \,\big] + (n-9) \delta _{9} + ... + (n-1) \delta _{1} .  $$
The subscripts of $ \delta  $ in the sum $ [\cdots ] $ are now consecutive.
Fortunately, this sum telescopes:
\begin{eqnarray*}
     \delta _{n-m-1} + ... + \delta _{10} &=& (-a_{m}+2a_{m+1}-a_{m+2}) + ... + (-a_{n-11}+2a_{n-10}-a_{n-9})      \\
     &=& -a_m+a_{m+1}+a_{n-10}-a_{n-9} \\
     &\leq & a_{n-10}-a_{n-9} .
\end{eqnarray*}
Hence,
\begin{equation}  \frac{2|T^-|}{x} \leq   \frac{(n-11)}{2} \, (a_{n-10}-a_{n-9})  + (n-9) \delta _{9} + ... + (n-1) \delta _{1} .   \Label{fin1}  \end{equation}
The RHS, using (\ref{gh9}) and (\ref{gh3}), adds up to less than
\begin{equation}  \frac{0.1636}{2} + 0.006902 + 0.010342 + 0.018 + 0.0326 + 0.3428 = 0.492444 .  \Label{fin2}  \end{equation}
This together with (\ref{fin1}), (\ref{fin2}), and (\ref{H4}) implies that $ S^- $ is NN, as
desired.

The proof of Theorem~AK is thus complete.

\end{enumerate}
}

\section{Lemmas}

\begin{LEMMA}[Fej\'er {[2, Satz XXVII]}] \Label{Fej}  {\ }

\par\vspace*{-3mm}\par
\begin{itemize}
\item[\rm(i)] If $ \left\{ c_1,c_2,\cdots ,c_m,c_{m+1}\right\}  $ is convex,
then $ \left[ c_1,\cdots ,c_m,\frac{c_{m+1}}{2} \right]  $ is NN in $ [0,\pi ] $.
\item[\rm(ii)] In particular, if $ \left\{ c_1,c_2,\cdots ,c_m,0\right\}  $ is convex,
then $ \left[ c_1,\cdots ,c_m\right]  $ is NN in $ [0,\pi ] $.
\end{itemize}
\end{LEMMA}

\begin{REM} \Label{cov2} \em
In (ii), assuming only the convexity of $ \left\{ c_1,c_2,\cdots ,c_m\right\}  $
is not enough to guarantee NN.
We need to augment 
the coefficient sequence with a 0 at the end. This is equivalent to requiring the
convexity of $ \left\{ c_1,c_2,\cdots ,c_m\right\}  $ plus $ a_{m-1}\geq 2a_m $.
\end{REM}

\begin{REM} \Label{cov} \em
It follows from the convexity of the power function $ x^\alpha  $
($ \alpha >1 $, $ x\geq 0 $) that if $ \{c_k^\beta \} $, $ \beta >0 $, is a convex sequence, so is $ \left\{ c_k^\gamma \right\}  $,
for $ \gamma >\beta  $.
\end{REM}

\begin{LEMMA}\Label{T3}
\begin{itemize}
\item[\rm(i)] For any integer $ k>0 $, $ \tau _{k-1}+\tau _{k} $ is NN in $ [0,\pi ] $.
\item[\rm(ii)] Let $ 0<A<B $. Then
$$  A\tau _{2j-1}^{-} + B\tau _{2j}^{-} \geq  -  j(B-A)  x \qquad  (x\geq 0)\,.  $$
\item[\rm(iii)] Let $ 0<m<m^* $ be integers, 
and $ d_k>0\ic $, $ k=m,\cdots ,m^* $. Then
$$  \sum_{k=m}^{m^*} d_k\tau _k \geq  d_{m^*} \min\left\{  \tau _{m^*}, 0 \right\}  .  $$
\item[\rm(iv)] In addition to the hypotheses of {\rm(iii)}, assume that
$ m $ is odd and $ m^* $ even. Then
\begin{equation}  \sum_{k=m}^{m^*} d_k\tau _k^{-} \geq  - \left[  \sum_{j=(m+1)/2}^{m^*/2} {j}\,(d_{2j}-d_{2j-1}) \right] \,x\,.  \Label{d2}  \end{equation}
\end{itemize}
\end{LEMMA}

\par\vspace*{2mm}\par
\begin{PROOF}
(i) This is a corollary of Lemma~\ref{Fej} (i), when $ c_k=2 $, $ k=1,\cdots ,m+1 $.
\begin{itemize}
\item[(ii)] By (\ref{tau2}),
$$  \tau _{2j}^- \geq   \frac{s\big(\frac{1}2\big)-s\big(2j+\frac12\big)}{2c\big(\frac12\big)} = \frac{-c\big(j+\frac{1}2\big)s\big(j\big)}{c\big(\frac12\big)}\,.  $$
In $ [0,\frac{\pi }{2j} ] $, $ c\big(j+\frac12\big)/c\big(\frac12\big) $ is a decreasing function.
Hence, it is less than its value at $ x=0 $, which is 1. It follows that 
$ \tau _{2j}^- \geq -s(j)\geq -jx $. In $ [\frac{\pi }{2j} ,\pi ] $, we obtain a lower bound of 
$ \tau _{2j}^- $ by replacing $ s\big(2j+\frac12\big) $ in the middle expression
above by $ -1 $ to obtain
$$  \tau _{2j}^- \geq   \frac{s\big(\frac{1}2\big)-1}{2c\big(\frac12\big)} \,.  $$
The RHS is negative but it is an increasing function of $ x $. At $ x=\frac{\pi }{2j} $, the RHS is greater
than $ -jx $, implying that the same is true for all $ x\in[\frac{\pi }{2j} ,\pi ] $. Thus
we have shown that $ \tau _{2j}^- \geq -jx $ for all $ x\in[0,\pi ] $.

By (i),
\begin{eqnarray*}
      A\tau _{2j-1}^- + B\tau _{2j}^- &=& A(\tau _{2j-1}^- +\tau _{2j}^-) + (B-A)\tau _{2j}^- \\
      &\geq & (B-A)\tau _{2j}^- \\
      &\geq & -j(B-A)x. 
\end{eqnarray*}
\item[(iii)] Rearrange the sum in question, in reverse order of the
terms, as
$$  d_{m^*} \tau _{m^*} + d_{m^*-1}\tau _{m^*-1} + ... + d_m\tau _m  $$

We compare this with the following sum, 
in which the coefficients are all the same:
$$  d_{m^*} ( \tau _{m^*} + \tau _{m^*-1} + ... + \tau _m ) .  $$
Let us look at the partial sums of the latter.
Using (i), we see that if the partial sum has an even number of terms, then it is NN. 
If there are an odd number of terms, the sum from the second term on is NN,
implying that the whole sum is not less than the first
term. In all cases, the partial sums 
$ {}\geq \min\left\{ d_{m^*}\tau _{m^*},0\right\}  $.
Applying the Comparison Principle then yields the conclusion.
The Comparison Principle is a well-known result. An explanation can be found in
\cite[See Lemma 2]{k2}.
\item[(iv)] This follows from (ii) and (iii).
\end{itemize}
\par\vspace*{-10mm}\par
{\ }\hfill
\end{PROOF}

\begin{LEMMA}\Label{ck}Suppose $ c_k>0,\dc $. Then
\begin{equation}  [c_1,c_2,\cdots ,c_n]^- \geq  \frac{c_1 \left[ s\big(\frac{1}2\big)+s\big(\frac{3}2\big)\right] -c_2\left[ 1+s\big(\frac{3}2\big)\right] }{2c\big(\frac{1}2\big)} \,.  \Label{cc1}  \end{equation}
\end{LEMMA}  \begin{PROOF}
The special case $ c_1=c_2=\cdots =1 $ corresponds to $ \tau _n^- $.  By (\ref{tau2}), its
partial sums satisfy
$$  \tau _m^-  \geq  \frac{s\big(\frac12\big)-1}{2c\big(\frac12\big)} \qquad  (1<m<n)\,.  $$
Applying the Comparison principle yields
\par\vspace*{-5mm}\par
$$  [1,1,c_2,\cdots ]^- \geq  \frac{s\big(\frac12\big)-1}{2c\big(\frac12\big)} \,,  $$
which is (\ref{cc1}) when $ c_1=c_2 $.
The general case follows from the relation
$$  [c_1,c_2,\cdots ,c_n]^- = (c_1-c_2)s(1)+[c_2,c_2,\cdots ,c_n]^-.  $$
\par\vspace*{-8mm}\par {\ }\hfill
\end{PROOF}

\begin{LEMMA}\Label{b5}
For $ x\in[0.75,\pi ] $,
$$  s\left( \frac{1}2\right) +0.338\,s\left( \frac{3}2\right) - 0.662 \geq  0  $$
\end{LEMMA}  \begin{PROOF}
Using the transformation $ x=\pi -2t $, we see that the desired inequality is equivalent to
$$  \cos(t)-0.338\cos(3t)-0.662 \geq  0, \qquad  \textstyle  t\in\big[0,\frac{\pi -0.75}{2} \big].  $$
With a further substitution $ X=\cos(t) $, the inequality becomes
$$  \frac{1}{500} (1-X)(676X^2+676X-331) \geq  0, \qquad  \textstyle  X \in\big[\cos\big(\frac{\pi -0.75}{2} \big),1\big] \subset [0.4,1],  $$
which is true, since all the roots of the LHS lie outside $ [0.366,1] $.
\end{PROOF}

Monotonicity properties of a sequence can often be deduced from corresponding
properties of its continuous analog. For instance, if $ \phi (x) $ is a differentiable
function in $ (0,\infty ) $, then $ \phi (k) $, $ k=1,2,\cdots  $ is a decreasing sequence if $ \phi '(x)\leq 0 $. Likewise,
\begin{eqnarray}
  \phi '(x)\ic\,(\dc)  &\Longrightarrow & \phi (k+1)-\phi (k)\ic\,(\dc) \nonumber  \\
  \phi ''(x)\ic\,(\dc)  &\Longrightarrow & \phi (k-1)-2\phi (k)+\phi (k+1)\ic\,(\dc) .  \Label{fpp}
\end{eqnarray}

\begin{LEMMA}\Label{fc}
{\rm(i)} For $ \beta \in[0,1] $.
$ \displaystyle \Box_k \dc \tx{in} k . $
\begin{itemize}
\item[\rm(ii)] Either $ \left\{ a_k\right\} _{k=1}^n $ is convex, or there is an $ m<n $
such that the sub-sequence
$ \left\{ a_1,\cdots ,a_m\right\}  $
is convex and the sub-sequence
$ \left\{ a_{m},a_{m+1},\cdots ,a_n\right\}  $
has an odd number of 
terms and is concave.
\item[\rm(iii)] When $ n\geq 7 $ and $ \beta \geq \beta _1 $, $ \Box_2,\Box_3,\Box_4\geq 0 $.
\end{itemize}

\end{LEMMA}  \begin{PROOF} 
(i) Note that
\begin{itemize}
\par\vspace*{-5mm}\par
\item[]
\begin{equation}  a_k = \left( \frac{n}{n^2-1} \right) ^\beta  f\left( \frac{k}{n} \right)  ,  \Label{akf}  \end{equation}
where
$$  \hspace*{6mm}        f(x) := \left( \frac{1}{x} - x \right) ^\beta  , \qquad  x \in [0,1].  $$
In view of (\ref{fpp}), the monotonicity of $ \Box_k $ will follow 
if we can show that for a fixed
$ \beta \in[0,1] $, $ f'''(x)<0 $ for $ x\in(0,1] $.
Direct computation gives
$$  f'''(x) = - \left( \frac{f(x)}{x^3(1-x^2)^3} \right)  g(x)  $$
where
{\small 
$$  g(x;\beta ) = \left( {x}^{6}+3{x}^{4}+3{x}^{2}+1 \right) {\beta }^{2}- \left( 3{x }^{6}+15{x}^{4}+9{x}^{2}-3 \right) \beta + \left( 2{x}^{6}+18{x}^{4}-6{x} ^{2}+2 \right)  .  $$
}

\par\vspace*{-5mm}\par
The conclusion follows if we can show that $ g(x)>0 $ for $ (x,\beta )\in[0,1]\times [0,1] $.
For each fixed $ x\in[0,1] $, $ g(x,\beta ) $, when extended to $ \beta \in(-\infty ,\infty ) $,
attains its minimum at
$$  \sigma  = \frac{3(x^4+4x^2-1)}{2(x^4+2x^2+1)} \,,   $$
which falls inside $ [0,1] $ only when $ x\in[x_1,x_2] $, where
$ x_1 = \sqrt{\sqrt5-2},\,\, x_2 = \sqrt{\sqrt{21}-4}. $

\underline{$\vphantom{y}$Case 1}: $ x\in[0,x_1] $. Minimum of $ g(x,\beta ) $ for $ \beta \in[0,1] $ is
attained when $ \beta =0 $. Hence,
$$  g(x;\beta ) \geq  g(x;0) = 2{x}^{6}+18{x}^{4}-6{x} ^{2}+2 > 0 .  $$
The positivity of the polynomial is verified using Sturm's
procedure.

\par\vspace*{3mm}\par
\underline{$\vphantom{y}$Case 2}: $ x\in[x_2,1] $. Minimum of $ g(x,\beta ) $ for $ \beta \in[0,1] $ is
attained when $ \beta =1 $. Hence,
$$  g(x;\beta ) \geq  g(x;1) = 6{x}^{4}-12{x}^{2}+6 > 0 .  $$

\underline{$\vphantom{y}$Case 3}: $ x\in(x_1,x_2) $. Minimum of $ g(x,\beta ) $ for $ \beta \in[0,1] $ is
attained when $ \beta =\sigma  $. Hence,
$$  g(x;\beta ) \geq  g(x;\sigma ) =   \frac{-x^8+8x^6-78x^4+56x^2-1}{4(x^2+1)} > 0 \,.   $$
Again positivity (for $x\in[0.4,0.8]\supset[x_1,x_2]$) is checked with
the Sturm procedure.

In all three cases, $ g(x;\beta )\geq 0 $; the first assertion 
of the Lemma is proved.

\item[(ii)] It is easy to verify $ \Box_2\geq 0 $.
If it happens that $ \Box_{k}\geq 0 $ for all $ k=2,\cdots ,n-1 $, then the sequence is 
convex. This situation prevails when $ \beta  $ is very close to 1.

In general, when $ n $ is large and $ \beta  $ is close to $ \beta _1 $,
$ f(x) $ is not convex in $ [0,1] $, as exemplified
by the red curve in Figure \ref{ff}.
Recall that after an appropriate scaling, see (\ref{akf}),
$ \hat a_n = f\big(\frac{k}{n} \big) $
are discrete points on the graph of $ f(x) $.
The curve starts out being convex and becomes concave after the point
of inflection, marked as a blue dot on the curve, between $ \hat{a}_{12} $ and
$ \hat{a}_{13} $. The two sub-sequences stipulated in the Lemma
are $ \left\{ a_1,\cdots ,a_{12}\right\}  $ and
$ \left\{ a_{13},\cdots ,a_{16},0\right\}  $, with the latter having five terms.

\begin{figure}[h]
\begin{center}
\fg{90}{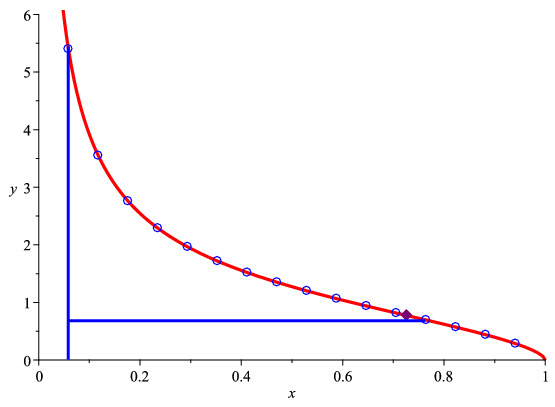}
\end{center}
\par\vspace*{-66.8mm}\par  \hspace*{42.5mm} $ \hat{a}_1 $
\par\vspace*{12mm}\par     \hspace*{47.0mm} $ \hat{a}_2 $
\par\vspace*{18.5mm}\par   \hspace*{92.4mm} $ \hat{a}_{12} $
\par\vspace*{0.9mm}\par   \hspace*{113.0mm} $ \hat{a}_{16} $

\par\vspace*{10mm}\par
\caption{\label{ff} Plot of $ f(x) $ (red curve) and $ \left\{ \hat{a}_{k}\right\}  $, $ n=17 $, $ \beta =\beta _1 $.}
\end{figure}

By (i), $ \Box $ will transition from $ \geq \,0 $ to $ <0 $ somewhere, as depicted
below:

$$  a_1, \,\, \cdots , \,\, \overbrace{\vp{7}{2} a_{\kappa  -1}, \, a_{\kappa } , \, a_{\kappa  +1} }^{\Box\geq 0}, \hspace*{-14mm} \underbrace{\vp{5}{5} \hspace*{15mm} a_{\kappa  +2} }_{\Box<0}, \,\, \cdots , \,\, a_n(=0)  $$

If we split the sequence right before $ a_{\kappa } $
or right after it, or after the next one, we get, in all three cases,
a convex sub-sequence
followed by a concave one. We can always choose one of these cases 
so that the second sub-sequence (including $ a_n=0 $)
has an odd number of terms. Letting $ a_m $ be the
last term of the first sub-sequence, we construct the two polynomials
represented, respectively, by the RHS of (\ref{SHK}) and (\ref{ST}). 

\item[(iii)] The monotonicity of $ \Box_k $ implies that, for a fixed $ \beta  $, there is a unique 
$ x_* $ at which $ f''(x_*)=0 $, and $ f''(x)>0 $ for $ x\in(0,x_*) $. Direct computation gives
$$  f''(x) = - \left( \frac{\beta f(x)}{x^2(1-x^2)^2} \right)  \left( (\beta -1)x^4+(2\beta -4)x^2+1\right)   ,  $$
and
$$  x_* = \sqrt{(\sqrt{5-4\beta }+\beta -2)/(1-\beta )}.  $$
As a function of $ \beta  $, $ x_* $ is increasing in $ \beta  $. Thus $ x_*\geq x_*(\beta _1)=0.5281747\cdots  $.
In terms of $ a_k $, this means that $ \Box_k\geq 0 $ for $ k=2,\cdots ,\kappa  $, where $ \kappa  $
is the largest integer less than or equal to $ 0.5281747n-1 $.
For $ n\geq 10 $, $ \kappa \geq 4 $ and the second assertion of the Lemma is proved. For
$ 5\leq n\leq 9 $, the assertion can be verified directly by computation.
\end{itemize}
\end{PROOF}

\par\vspace*{1mm}\par
\begin{LEMMA}\Label{ABC}
Let $ \lambda ,\mu \geq 0 $, and $ 0<C<B<A\leq 1 $. The function
\begin{equation}  \xi (\beta )=\lambda A^\beta -(\lambda +\mu )B^\beta +\mu C   \Label{xi}  \end{equation}
is $ \geq \,0 $ at some point $ \beta _0>0\,\Longrightarrow \,\xi (\beta ) $ cannot have a 
local minimum in $ (\beta _0,\infty ) $.
\end{LEMMA}
\par\vspace*{2mm}\par
\begin{PROOF}
$ \xi (\beta _0)\geq 0 $ $ \Longrightarrow  $ $ \xi (\beta )>0 $ for $ \beta >\beta _0 $
(due to the convexity of $ x^\beta  $, $ \beta >1 $).

A local minimum must be a critical point.
Our tactic is to show that at a critical point $ \sigma  $ (i.e. when
$\xi '(\sigma )=0$), $ \xi ''(\sigma )<0 $,
which implies that the critical point cannot be a local minimum.

We can assume, without loss of generality, that $ \sigma =1 $
(use $ A^{\sigma } $ and $ B^{\sigma } $ as the new $ A $
and $ B $, respectively). That $ 1 $ is a critical point gives
\begin{equation}  \xi '(1) = \lambda  A\ln(A)-(\lambda +\mu )B\ln(B) + C\ln(C) = 0  \Label{alna0}  \end{equation}
and we need to show
\begin{equation}  \xi ''(1) = \lambda  A\ln^2(A)-(\lambda +\mu )B\ln^2(B) + C\ln^2(C) < 0 .  \Label{alna1}  \end{equation}

The point $ \sigma =\mbox{e}^{-1}\approx 0.3678794412 $ divides $ [0,1] $ into
two sub-intervals: the function $ f(t)=|t\ln(t)| $ is $ \ic $ in $ [0,\sigma ] $,
but $ \dc $ in $ [\sigma ,1] $.
In each sub-interval, $ f(t) $ is invertible. Let $ g(s)=f^{-1}(s) $ be
the inverse of $ f(t) $ in $ [0,\sigma ] $, i.e. $ |g(s)\ln(g(s))|=s $.

We divide the proof of (\ref{alna1}) into three cases.

\underline{$\vphantom{y}$Case 1}: $ A\in[0,\sigma ] $. Then $ B $ and $ C $ are also in $ [0,\sigma ] $.

That (\ref{alna1}) follows from (\ref{alna0}) is a consequence of the concavity of 
$ t\ln^2(t)=g(s)\ln^2(g(s)) $ as a function of $ s=t\ln(t) $. Using the
chain rule, one can verify that
$$  \frac{d}{ds} (t\ln^2(t)) = \frac{d}{dt} (t\ln^2(t)) \frac{dt}{ds} = \frac{\ln^2(t)+2\ln(t)}{\ln(t)+1}   $$
and
$$  \frac{d^2}{ds^2} (t\ln^2(t)) = \frac{d}{dt} \left(  \frac{\ln^2(t)+2\ln(t)}{\ln(t)+1} \right)  \frac{dt}{ds} = \frac{\ln^2(t)+2\ln(t)+2}{t(\ln(t)+1)} \,.  $$
The last expression is $ - $  for $ t\in[0,\sigma ] $
(the numerator is +, but the denominator is $ - $), proving
the claim.

\begin{figure}[h]
\begin{center}
\fg{70}{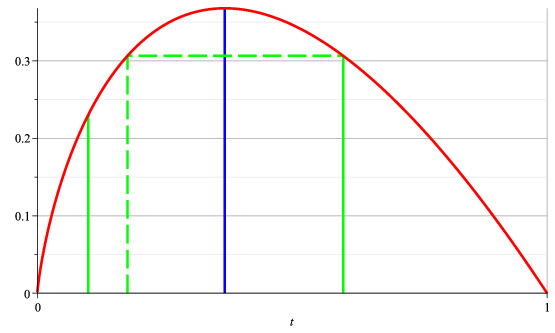}
\end{center}
\par\vspace*{-10.5mm}\par  
\hspace*{37mm} $ C $\hspace*{3mm} $ A_1 $\hspace*{33.5mm} $ A $
\par\vspace*{1mm}\par
\caption{\label{ff2} Graph of $ f(t)=|t\ln(t)| $. Case 1, sub-case 1.}
\end{figure}

\underline{$\vphantom{y}$Case 2}: $ A\in[\sigma ,1] $ and $ C\in[0,\sigma ] $.
There are two sub-cases. The first is $ f(A)\geq f(C) $. One such 
example is depicted in Figure \ref{ff2}. There exists $ A_1\in[C,\sigma ] $ such that
$ f(A_1)=f(A) $. Then, with $ A_1 $ replacing $ A $, we have the same
situation as Case 1. Thus, (\ref{alna1}), with $ A_1 $ in place of $ A $, i.e.
\begin{equation}  \lambda  A_1\ln^2(A_1)-(\lambda +\mu )B\ln^2(B) + C\ln^2(C) < 0 .  \Label{alna3}  \end{equation}
Since $ \ln(A_1)<\ln(A)<0 $,
\begin{equation}  A\ln^2(A) = (A\ln(A))\ln(A) =(A_1\ln(A_1)\ln(A) < (A_1\ln(A_1))\ln(A_1).  \Label{alna4}  \end{equation}
Now (\ref{alna3}) and (\ref{alna4}) imply (\ref{alna1}).

\begin{figure}[h]
\begin{center}
\fg{70}{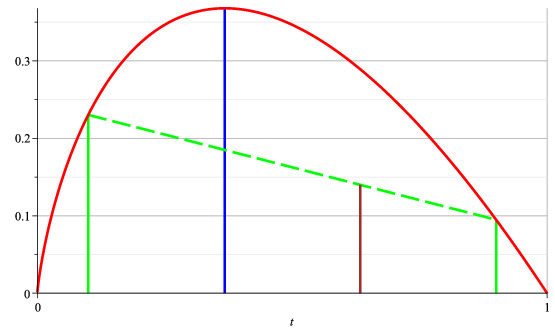}
\end{center}
\par\vspace*{-10.5mm}\par  
\hspace*{37mm} $ C $\hspace*{45mm} $ D $\hspace*{20.0mm} $ A $
\par\vspace*{-31mm}\par
\hspace*{86mm} $ E $
\par\vspace*{28mm}\par
\caption{\label{ff3} Graph of $ f(t)=|t\ln(t)| $. Case 1, sub-case 2.}
\end{figure}

The second sub-case is when $ f(A)<f(C) $, as depicted
in Figure \ref{ff3}. Let $ D=(\lambda A+\mu C)/(\lambda +\mu ) $; it lies between $ C $ and $ A $.
The assumption $ \xi (1)\geq 0\,\Longrightarrow \,B\in[C,D] $.

The dash green line joins $ (C,f(C)) $ and $ (A,f(A)) $.
Concavity of $ f(t)\,\Longrightarrow \,{} $ dash green line $ \leq  $ red curve.
Let the vertical line through $ D $ cut the dash green line at $ E $. The
length $ DE $ is $ (\lambda f(A)+\mu f(C))/(\lambda +\mu )=f(B) $, by (\ref{alna0}).
However, this contradicts the earlier assertion that $ B $ lies between $ C $ and
$ D $ because for all $ B\in[C,D] $, $ B\ln(B) $ (the red curve) is obviously
larger than $ DE $ (the red curve is above $ E) $. This contradiction means
that (\ref{alna0}) cannot hold. In other words, $ 1 $ cannot be a critical point
of $ \xi (b) $.

\underline{$\vphantom{y}$Case 3}: $ C\in[\sigma ,1] $.
Then $ f(A)<f(C)\,\Longrightarrow \,{} $ contradiction just as in the last case.
\end{PROOF}

A simple corollary is

\par\vspace*{1mm}\par
\begin{LEMMA}\Label{AlnA2}Let $ \xi (\beta ) $ be defined as in Lemma~\ref{ABC}
and $ \xi (\beta _1)\geq 0 $.
\begin{itemize}
\item[\rm(i)] If $ \xi '(1)\geq 0 $, then $ \xi (\beta ) $ $ \ic $ in $ [\beta _1,1] $.
\item[\rm(ii)] If $ \xi '(\beta _1)\leq 0 $, then $ \xi (\beta ) $ $ \dc $ in $ [\beta _1,1] $.
\end{itemize}
\end{LEMMA}

\par\vspace*{2mm}\par
\begin{PROOF}
(i) Suppose the contrary. 
Then $ \xi '(\beta ) $ must be negative at some point $ \beta _3\in[\beta _1,1) $, and
\par\vspace*{-5.5mm}\par
\begin{itemize}
\item[{}]
there must exist a local minimum between $ \beta _3 $ and 1, contradicting the Lemma.

\item[(ii)] In a similar way, suppose the contrary.
Then $ \xi '(\beta ) $ must be positive at some point $ \beta _3\in[\beta _1,1) $, and
there must exist a local minimum between $ \beta _3 $ and 1, a contradiction.
\end{itemize}
\par\vspace*{-5mm}\par \hfill
\end{PROOF}

The coefficients of $ H $, as defined in (\ref{SH}), written in terms of 
$ y\in\big[0,\frac{1}{48} \big] $ and $ \beta \in[\beta _1,1] $ are given by $ h_i(y;\beta ) $,
$ i=1,2,3 $ below, while the second differences are given by $ i=4,5 $.
\begin{eqnarray*}
    h_1(y;\beta ) &:=& 1 - 4 \left( \frac{1-15y}{4} \right) ^\beta  + 3 \left( \frac{1-24y}{5} \right) ^\beta  \\
    h_2(y;\beta ) &:=& \left( \frac{1-3y}{2} \right) ^\beta  - 3 \left( \frac{1-15y}{4} \right) ^\beta  + 2 \left( \frac{1-24y}{5} \right) ^\beta  \\
    h_3(y;\beta ) &:=& \left( \frac{1-8y}{3} \right) ^\beta  - 2 \left( \frac{1-15y}{4} \right) ^\beta  +  \left( \frac{1-24y}{5} \right) ^\beta   \\
h_4(y;\beta ) &:=& h_1(y;\beta )-2h_2(y;\beta )+h_3(y;\beta ) \\
          &=& 1 - 2 \left( \frac{1-3y}{2} \right) ^\beta  +  \left( \frac{1-8y}{3} \right) ^\beta  \\
h_5(y;\beta ) &:=& h_2(y;\beta )-2h_3(y;\beta ) \\
          &=&  \left( \frac{1-3y}{2} \right) ^\beta  - 2 \left( \frac{1-8y}{3} \right) ^\beta  +  \left( \frac{1-15y}{4} \right) ^\beta . \hspace*{4.5mm}   
\end{eqnarray*}

\begin{LEMMA}\Label{hhh}
\begin{itemize}
\item[\rm(i)] With $ \beta =\beta _1 $, all five functions $ h_i(y;\beta _1) $ $(i=1,\cdots ,5)$, 
are + $ \dc $ in $ y $.
\item[\rm(ii)] For fixed $ y $, all five functions
$ h_i(y;\beta ) $ $(i=1,\cdots ,5)$, 
are $ \ic $ in $ \beta  $.
\end{itemize}
\end{LEMMA}

\begin{PROOF}
(i) All conclusions are proved using the MAPLE procedure {\tt dif}, see \cite{k1}.

\noindent \underline{$\vphantom{y}$$h_1(y;\beta _1)$}: \quad We want to show that
\begin{eqnarray*}
   - \, \frac{5\,h_1'(y;\beta _1)}{\beta _1} &=&   72   \left( \frac{1-24y}{5} \right) ^{\beta _1-1}  - 75   \left( \frac{1-15y}{4} \right) ^{\beta _1-1} 
\end{eqnarray*}
is positive for $ y\in\big[0,\frac{1}{48} \big] $.
The MAPLE display is shown in Figure \ref{fb1}.

\newpage

\par\vspace*{-15mm}\par
\begin{center}
\fbox{\begin{minipage}{5.4in}\vspace*{5pt}\includegraphics[width=140mm]{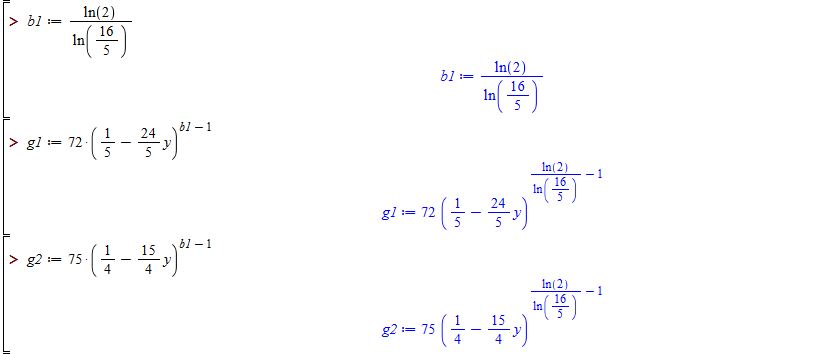}\par\includegraphics[width=140mm]{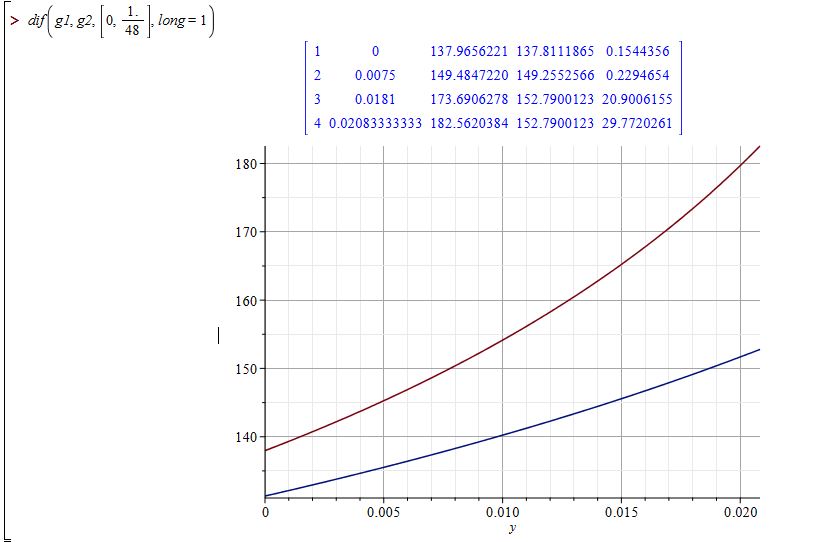}\vspace*{5pt}\end{minipage}}

\par\vspace*{5mm}\par
Figure \F. \Label{fh1} MAPLE output for $ h_1(y;\beta _1) $.
\end{center}

The first command defines {\tt b1} ($\beta _1$). 
The next two commands defines {\tt g1} and {\tt g2}. The fourth
invokes {\tt dif} in the verbose format (by adding the option {\tt long=1})
to generate four points (the second column
of the displayed matrix):
$$  \tau _1=0,\,\, \tau _2 = 0.0075,\,\, \tau _3=0.0181, \,\, \tau _4=\textstyle  \frac{1}{48} .  $$
The third column of the matrix lists {\tt g1}$(\tau _i)$, $(i=1,2,3,4)$
the fourth lists {\tt g2}$(\tau _i)$, $(i=2,3,4,4)$, while the fifth column
is the third column minus the fourth and must be positive. 
The successful generation of the matrix proves that {\tt g1 - g2} is
positive in the interval under study. The graphs are plot as a visual check that
{\tt g1} (the red curve) and {\tt g2} (the blue curve) are monotone.

\noindent \underline{$\vphantom{y}$$h_2(y;\beta _1)$}: \quad We have
$$  - \, \frac{h_2'(y;\beta _1)}{\beta _1} = \left[  \frac{3}{2}  \left( \frac{1-3y}{2} \right) ^{\beta _1-1} + \frac{48}{5} \left( \frac{1-24y}{5} \right) ^{\beta _1-1} \right]  - \frac{45}{4} \left( \frac{1-15y}{4} \right) ^{\beta _1-1} \,.  $$
This time, we invoke {\tt dif} without the {\tt long=1} option. The output
then only shows the list of $ \tau _i $, not the verbose matrix or the graphs.

\par\vspace*{8mm}\par
\begin{center}
\FG{140}{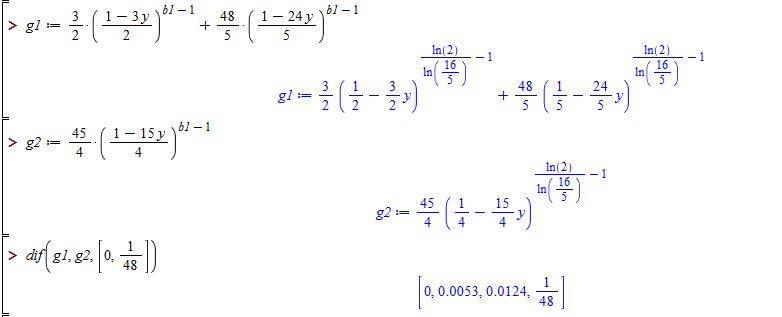} 

\par\vspace*{5mm}\par
Figure \F. \Label{fh2} MAPLE output for $ h_2(y;\beta _1) $.
\end{center}

\par\vspace*{5mm}\par
\noindent \underline{$\vphantom{y}$$h_3(y;\beta _1)$}: \hspace*{5mm}
$ \displaystyle - \, \frac{h_3'(y;\beta _1)}{\beta _1} = \left[  \frac{8}{3}  \left( \frac{1-8y}{3} \right) ^{\beta _1-1} + \frac{24}{5} \left( \frac{1-24y}{5} \right) ^{\beta _1-1} \right]  - \frac{15}{2} \left( \frac{1-15y}{4} \right) ^{\beta _1-1} \,. $

\begin{center}
\FG{140}{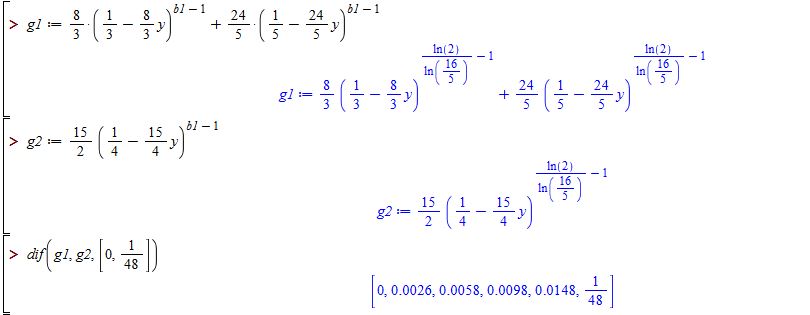} 

\par\vspace*{5mm}\par
Figure \F. \Label{fh3} MAPLE output for $ h_3(y;\beta _1) $.
\end{center}

\noindent \underline{$\vphantom{y}$$h_4(y;\beta _1)$}: 

\begin{center}
\FG{140}{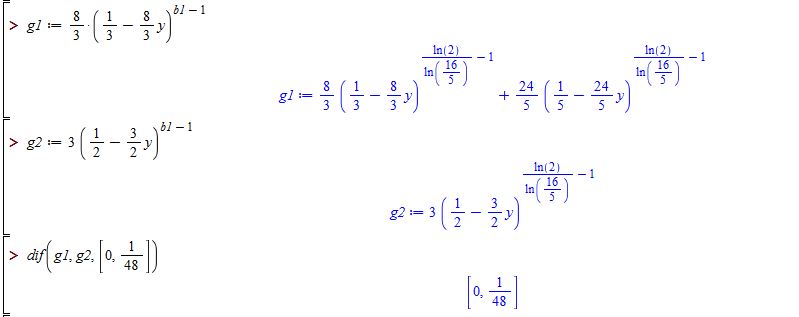} 

\par\vspace*{5mm}\par
Figure \F. \Label{fh4} MAPLE output for $ h_4(y;\beta _1) $.
\end{center}

\noindent \underline{$\vphantom{y}$$h_5(y;\beta _1)$}: 

\begin{center}
\FG{140}{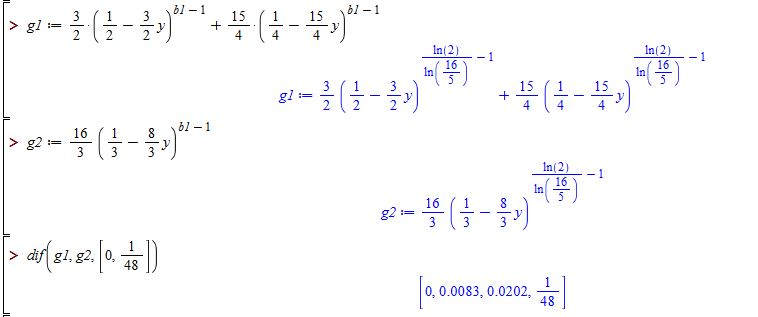} 

\par\vspace*{5mm}\par
Figure \F. \Label{fh5} MAPLE output for $ h_5(y;\beta _1) $.
\end{center}

\noindent (ii) Although we assume that $ y $ is fixed, it is not given a 
specific value. In other words, it is used as a parameter,
and for that reason, {\tt dif} cannot be applied directly.
Instead, we make use of Lemma~\ref{AlnA2}, noting that $ h_i(y;\beta ) $
(with fixed $ y $) has the same form as $ \xi (\beta ) $ in (\ref{xi}).
By Lemma~\ref{AlnA2} (i), it suffices to show that 
$$  \frac{\partial }{\partial \beta } \, h_i(y;1) \geq  0, \quad  i=1, \cdots , 5.  $$
Now we can apply {\tt dif} to each $ \partial h_i(y;1)/\partial \beta  $.

\newpage
\par\vspace*{-10mm}\par
\noindent \underline{$\vphantom{y}$$\partial h_1(y;1)/\beta  $}: \hspace*{10mm}
$ \displaystyle 3\left( \frac{1-24y}{5} \right)  \ln \left( \frac{1-24y}{5} \right)  - 4\left( \frac{1-15}{4} \right)  \ln \left( \frac{1-15}{4} \right)  $

\begin{center}
\FG{140}{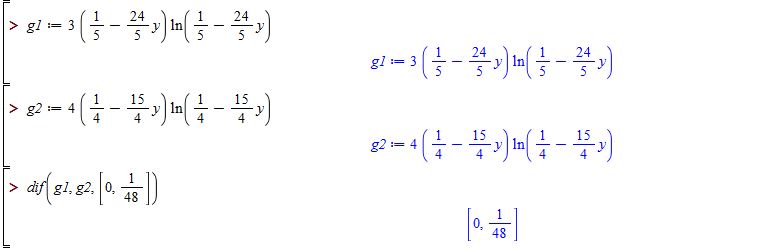} 

\par\vspace*{5mm}\par
Figure \F. \Label{fph1} MAPLE output for $ \partial h_1(y;1)/\partial \beta  $.
\end{center}

\par\vspace*{-3mm}\par
\noindent \underline{$\vphantom{y}$$\partial h_2(y;1)/\beta  $}:

\begin{center}
\FG{140}{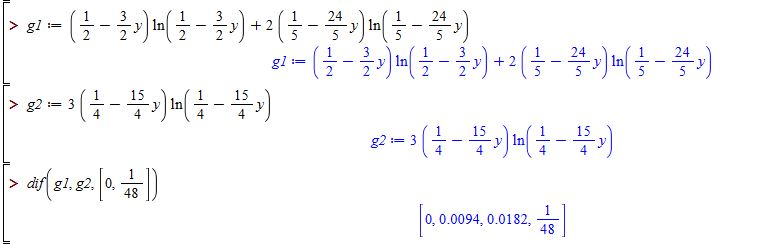} 

\par\vspace*{5mm}\par
Figure \F. \Label{fph2} MAPLE output for $ \partial h_2(y;1)/\partial \beta  $.
\end{center}

\par\vspace*{-3mm}\par
\noindent \underline{$\vphantom{y}$$\partial h_3(y;1)/\beta  $}:

\par\vspace*{-3mm}\par
\begin{center}
\FG{140}{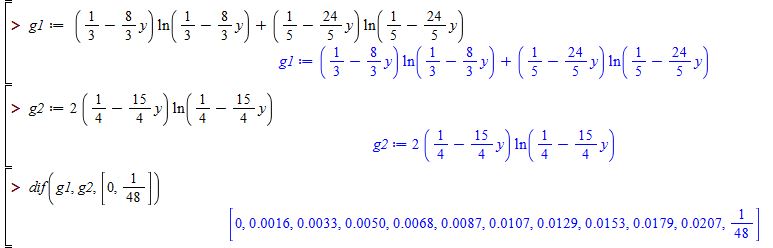} 

\par\vspace*{5mm}\par
Figure \F. \Label{fph3} MAPLE output for $ \partial h_3(y;1)/\partial \beta  $.
\end{center}

\noindent \underline{$\vphantom{y}$$\partial h_4(y;1)/\beta  $}:

\begin{center}
\FG{140}{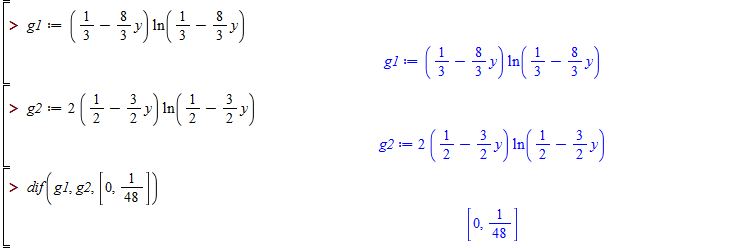} 

\par\vspace*{5mm}\par
Figure \F. \Label{fph4} MAPLE output for $ \partial h_4(y;1)/\partial \beta  $.
\end{center}

\noindent \underline{$\vphantom{y}$$\partial h_5(y;1)/\beta  $}:

\begin{center}
\FG{140}{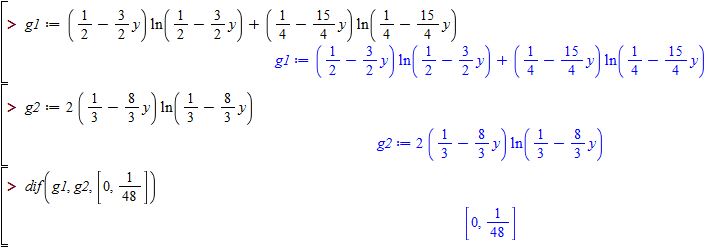} 

\par\vspace*{5mm}\par
Figure \F. \Label{fph5} MAPLE output for $ \partial h_5(y;1)/\partial \beta  $.
\end{center}

\par\vspace*{-10mm}\par
{ }\hfill
\end{PROOF}

\begin{LEMMA}\Label{BB}
For any $ 0<B<1 $, $ \beta \in\big[\frac{1}{2} ,1\big] $,
$ (1+B)^\beta  + (1-B)^\beta  \ic  $
in $ \beta  $.
\end{LEMMA}

\begin{PROOF}
The conclusion is true if we can show that
\begin{eqnarray}
   && \frac{\partial }{\partial \beta } \left[  (1+B)^\beta  + (1-B)^\beta  \right]  \nonumber  \\
   &=& (1+B)^\beta  \ln(1+B) + (1-B)^\beta  \ln(1-B) \nonumber  \\
   &\geq & 0 .   \Label{dB}
\end{eqnarray}
Since $ \beta >\frac{1}{2} $, (\ref{dB}) follows from the stronger assertion
$$  (1+B)^{1/2} \ln(1+B) + (1-B)^{1/2} \ln(1-B) \geq  0 ,  $$
which is equivalent to
$$  \theta _1(B)-\theta _2(B) := (1+B)\ln^2(1+B)-(1-B)\ln^2(1-B) \geq  0 .     $$
We divide the proof into three cases:
\begin{itemize}
\item[(1)] $ \theta _1(B) $ is increasing in $ B\in[0,1] $. $ \theta _2(B) $ is increasing in $ [0,\tau ] $
  and decreasing in $ [\tau ,1] $, where $ \tau =1-\mbox{e}^{-2} $. Since 
  $ \theta _1(\tau )>\theta _2(\tau ) $, we have $ \theta _1(B)>\theta _2(B) $
  in $ [\tau ,1] $.
\item[(2)] In $ [0.4,\tau ] $, we can use the DIF technique 
  to confirm that $ \theta _1(B)>\theta _2(B) $.
\item[(3)] The Taylor series of $ \ln(1+x) $ is an alternating (i.e. $ + $ and $ - $)
  series. Hence,
$$  \ln(1+B) \geq  B - \frac{B^2}{2} + \frac{B^3}{3} - \frac{B^4}{4}   $$
$\,\Longrightarrow \,$
\begin{equation}  \theta _1(B) \geq  (1+B) \left(  \frac{B^2}{2} + \frac{B^3}{3} - \frac{B^4}{4} \right) ^2 .  \Label{h1B}  \end{equation}
On the other hand, the Taylor series of $ \theta _2(B) $ is
$$  {B}^{2}-{\frac{1}{12}}{B}^{4}-{\frac{1}{12}}{B}^{5}-{\frac{13}{180}}{ B}^{6}-{\frac{11}{180}}{B}^{7}-{\frac{29}{560}}{B}^{8}-{\frac{223}{ 5040}}{B}^{9}+ ... ,  $$
suggesting that
\begin{equation}  {B}^{2}-{\frac{1}{12}}{B}^{4}-{\frac{1}{12}}{B}^{5} \geq  \theta _2(B) .  \Label{h2B}  \end{equation}
This is confirmed because, for $ B\in[0,0.4] $,
$$  \frac{d^6}{dB^6} \left( {B}^{2}-{\frac{1}{12}}{B}^{4}-{\frac{1}{12}}{B}^{5} - \theta _2(B) \right)  = \frac{52-48\ln(1-B)}{(1-B)^5} \geq  0 .  $$
         
The Sturm procedure $ \,\Longrightarrow \,{} $ for $ B\in[0,0.4] $,
\begin{equation}  (1+B) \left(  \frac{B^2}{2} + \frac{B^3}{3} - \frac{B^4}{4} \right) ^2 \geq   {B}^{2}-{\frac{1}{12}}{B}^{4}-{\frac{1}{12}}{B}^{5} .  \Label{h3B}  \end{equation}
(\ref{h1B}), (\ref{h2B}), + (\ref{h3B})  $ \,\Longrightarrow \,\theta _1(B)\geq \theta _2(B) $ in 
$ [0,0.4] $.
\end{itemize}
This completes the proof of the Lemma.
\end{PROOF}

\begin{LEMMA}\Label{H7b}
The function
\begin{eqnarray*}
     && \textstyle  \left( 1 - 4 \left( \frac{11}{64} \right) ^{\beta _1} + 3 \left( \frac{1}{10} \right) ^{\beta _1} \right) \sin(x) - \left(  \left( \frac{15}{32} \right) ^{\beta _1}- 3 \left( \frac{11}{64} \right) ^{\beta _1} + 2 \left( \frac{1}{10} \right) ^{\beta _1}  \right) \sin(2x) \\
     && \textstyle  \hspace*{3mm}  + \left( \left( \frac{5}{18} \right) ^{\beta _1} - 2 \left( \frac{11}{64} \right) ^{\beta _1} +  \left( \frac{1}{10} \right) ^{\beta _1} \right)  \sin(3x) 
\end{eqnarray*}
is concave in $ [0,0.75] $.
\end{LEMMA}  \begin{PROOF}
Denote the function by $ \varphi (x) $. We need to show that $ \varphi ''(x)<0 $ in the given interval.
In terms of $ X=\cos(x)\in[\cos(0.75),1] $,
\begin{eqnarray*}
      \frac{\varphi ''(x)}{\sin(x)} \! \! &=& \! \! \textstyle  \left( 72 \left( \frac{11}{64} \right) ^{\beta _1} - 36 \left( \frac{1}{10} \right) ^{\beta _1}  - 36 \left( \frac{5}{18} \right) ^{\beta _1} \right)  X^2 + \left(  8 \left( \frac{15}{32} \right) ^{\beta _1}- 24 \left( \frac{11}{64} \right) ^{\beta _1} + 16 \left( \frac{1}{10} \right) ^{\beta _1}  \right) X \\
     && \textstyle  \hspace*{3mm}  + \left(  9 \left( \frac{5}{18} \right) ^{\beta _1} - 14 \left( \frac{11}{64} \right) ^{\beta _1} + 6 \left( \frac{1}{10} \right) ^{\beta _1} -1 \right)  .
\end{eqnarray*}

\noindent It is straightforward to check that the roots of the RHS are
$$  0.39281956258689586, \quad  0.67755077339437549,  $$
computed using MAPLE with an accuracy of 40 digits, both of which
lie outside the given interval. Hence, $ \varphi ''(x) $ is of one sign in the interval,
and the sign is negative, as desired.
\end{PROOF}

\begin{LEMMA}
For $ n\geq 7 $,\quad 
$ \displaystyle n \left(  \frac{2n-1}{(n^2-1)(n-1)} \right) ^{\beta _1}  $
\quad  is a decreasing function of $ n $.
\end{LEMMA}  \begin{PROOF}
Since $ \beta _1>\frac{1}{2} $, the conclusion follows if we can show that
$$  n \left(  \frac{2n-1}{(n^2-1)(n-1)} \right) ^{1/2}   $$
is decreasing in $ n $, or equivalently, if
$$  \frac{n^2(2n-1)}{(n^2-1)(n-1)}  $$
is decreasing in $ n $. This is a routine exercise in calculus.
\end{PROOF}

 
Let $ \delta _k := - \Box_{n-k}$ be as defined in Section~2 \P12.
Then,
\begin{eqnarray*}
         \delta _{1} &=& \left( \frac{1}{n^2-1} \right) ^\beta  \left[  2\left( \frac{2n-1}{n-1} \right) ^\beta  - \left( \frac{4n-4}{n-2} \right) ^\beta  \right]    \\
          &=& \left( \frac{2n-1}{(n^2-1)(n-1)} \right) ^\beta   \left[  2 - \left( \frac{4(n-1)^2}{(2n-1)(n-2)} \right) ^\beta  \right]   . 
\end{eqnarray*}
\begin{LEMMA}
\begin{equation}  (n-1)\delta _1 \leq  \theta (12) < 0.3921 \qquad \tx{for} n\geq 7.  \Label{gh1}  \end{equation}
\begin{equation}  (n-1)\delta _1 \leq  \theta (45) < 0.3428 \qquad \tx{for} n\geq 45.  \Label{gh2}  \end{equation}
\end{LEMMA}

\begin{PROOF}
Since
$ \frac{2n-1}{(n^2-1)(n-1)} < 1 , $
the first factor $ \dc $ in $ \beta  $.
On the other hand,
\mbox{ $ \frac{(4n-4)(n-1)}{(2n-1)(n-2)} > 1 $, }
implying that the second factor is also $ \dc $ in $ \beta  $.
It follows that
\begin{equation}  (n-1) \delta _{1} \leq  (n-1) \left( \frac{2n-1}{(n^2-1)(n-1)} \right) ^{\beta _1}  \left[  2 - \left( \frac{4(n-1)^2}{(2n-1)(n-2)} \right) ^{\beta _1} \right]  .  \Label{dn2}  \end{equation}
Denote by $ \theta (n) $ the RHS of (\ref{dn2}).  Numerics gives
$$  \theta (7)<\theta (8)<\theta (9)<\theta (10)<\theta (11)<\fb{\theta (12)}>\theta (13)>\theta (14)> ...  $$

We rewrite
\begin{equation}  \theta (n) = \left(  \frac{2n-1}{(n+1)(n-1)^{0.1}} \right) ^{\beta _1} \raisebox{4pt}{$\left [ \vp{17}{0} \right . $} \frac{ 2 - \left( \frac{4(n-1)^2}{(2n-1)(n-2)} \right) ^{\beta _1} }{(n-1)^{1.9\beta _1-1}}  \raisebox{4pt}{$\left . \vp{17}{0} \right ] $} \,.  \Label{hn}  \end{equation}
Direct computation gives
$$  \frac{d}{dn} \left(  \frac{2n-1}{(n+1)(n-1)^{0.1}} \right)  = - \, \frac{2n^2-29n+29}{10(n+1)^2(n-)^{1.1}} < 0 \qquad  \tx{when} n > 14.  $$
It follows that the first factor in (\ref{hn}) is $ \dc $ ($ n>14 $).
Rewrite the second factor as
$ \displaystyle \frac{\eta (y)}{y^\alpha } ,  $
where
$ y=n-1, \eta (y) = 2 - \left( \frac{4y^2}{(2y+1)(y-1)} \right) ^{\beta _1}, \alpha =1.9\beta _1-1 . $
Obviously,
$ \eta (y) $ is $ \ic $ in $ y>12 $.

Direct computation gives
$$  y\eta '(y) = \beta _1 \left( \frac{y+2}{2y+1} \right)  \left(  \frac{4y^2}{(2y+1)(y-1)^{1+1/\beta _1}} \right) ^{\beta _1} .  $$
Both fractions are $ \dc $ in $ y\,\Longrightarrow \,y\eta '(y) $
$ \dc $ in $ y $.
$$  \left(    \frac{\eta (y)}{y^\alpha } \right) ' = \frac{y\eta '(y)-\alpha \eta (y)}{y^{\alpha +1}} \,.       $$
The numerator, being a $ \dc $ function minus an $ \ic $ function, is $ \dc $ in~$y$. 
Substituting
$ y=10 $ shows that the numerator is negative. Hence, the numerator is negative for all
larger $ y $. Thus the second factor in (\ref{hn}) is $ \dc $, implying that $ \theta (n) $
is $ \dc $ in $ n>14 $. Direct computation shows $ \theta (12)>\theta (13)>\theta (14) $. Thus 
$ \theta (n) $ is $ \dc $ for $ n\geq 12 $.
\end{PROOF}

\par\vspace*{2mm}\par
\begin{LEMMA}
For odd $ n\geq 15 $,

\par\vspace*{-6mm}\par
\begin{eqnarray}
   (n-3) \delta _{3}  &<&  0.0412     \Label{gh3}  \\
   (n-5) \delta _{5}  &<&  0.018      \Label{gh5}  \\
   (n-7) \delta _{7}  &<&  0.010342   \Label{gh7}  \\
   (n-9) \delta _{9}  &<&  0.006902 .  \Label{gh9} 
\end{eqnarray}

For $ n\geq 45 $
\begin{equation}  (n-3) \delta _3 \hspace*{2.5mm} < \hspace*{2.5mm} 0.0326 .  \hspace*{3mm}  \Label{dn6}  \end{equation}
\end{LEMMA}

\begin{PROOF}
For (\ref{gh3})--(\ref{gh9}), we only give the proof of (\ref{gh3}). 
The other inequalities can be established in exactly the same way.

\begin{equation}  \delta _{3} = \left( \frac{6n-9}{(n^2-1)(n-3)} \right) ^{\beta }  \left[  2 - \left( \frac{(8n-16)(n-3)}{(6n-9)(n-4)} \right) ^{\beta } - \left( \frac{(4n-4)(n-3)}{(6n-9)(n-2)} \right) ^{\beta } \right]  .  \Label{dn3}  \end{equation}
The first factor is $ \dc $ in $ \beta  $. We claim that the
second factor is also $ \dc $, which is equivalent to 
$$  \left( \frac{(8n-16)(n-3)}{(6n-9)(n-4)} \right) ^{\beta } + \left( \frac{(4n-4)(n-3)}{(6n-9)(n-2)} \right) ^{\beta } = (1+A)^\beta  + (1-B)^\beta   $$
being $ \ic $ in $ \beta  $. It is easy to verify that
$$  0 <    B = \frac{2n^2-5n+6}{3(2n-3)(n-2)} < \frac{2n^2-7n+12}{3(2n-3)(n-4)} = A < 1 .  $$
Obviously, 
$ (1+A)^\beta -(1+B)^\beta  \ic $
ic in $ \beta  $. The claim then follows from this and Lemma~\ref{BB}.

It then follows that
\begin{eqnarray}
\delta _{3} &\leq & \left( \frac{6n-9}{(n^2-1)(n-3)} \right) ^{\beta _1}  \left[  2 - \left( \frac{(8n-16)(n-3)}{(6n-9)(n-4)} \right) ^{\beta _1} - \left( \frac{(4n-4)(n-3)}{(6n-9)(n-2)} \right) ^{\beta _1} \right]  \nonumber \\
   &=& F(n) G(n) ,    \Label{dn4}
\end{eqnarray}
where $ F(n) $ and $ G(n) $ are the first and second factor of (\ref{dn3}), respectively. 
Obviously, $ F(n) $ is $ \dc $ in $ n $.

Our next claim is that $ G(n) $ is $ \ic $ in $ n\geq 15 $. 
This is equivalent to 
$$  \left( \frac{(8n-16)(n-3)}{(6n-9)(n-4)} \right) ^{\beta _1} + \left( \frac{(4n-4)(n-3)}{(6n-9)(n-2)} \right) ^{\beta _1} \, \, \dc .  $$
After differentiating this expression and canceling some common factors, 
we see that this is equivalent to
$$  - 2 \left( \frac{n^2-6}{(n-4)^2} \right) \left( \frac{8n-16}{n-4} \right) ^{\beta _1-1} + \left( \frac{n^2-3}{(n-2)^2} \right) \left( \frac{4n-4}{n-2} \right) ^{\beta _1-1} \leq  0   $$
which is, in turn, equivalent to
\begin{equation}  \left( \frac{(8n-16)(n-2)}{(4n-4)(n-4)} \right) ^{1-\beta _1} \leq  \frac{2(n^2-6)(n-2)^2}{(n^2-3)(n-4)^2} \,.  \Label{fg1}  \end{equation}
Since the fraction inside the parentheses on the LHS is $ {}>1 $, (\ref{fg1}) follows
from
$$  \frac{(8n-16)(n-2)}{(4n-4)(n-4)}  \leq  \frac{2(n^2-6)(n-2)^2}{(n^2-3)(n-4)^2} \,,  $$
which is easy to verify. The proof of the claim is complete.

Then (\ref{dn4}) leads to
\begin{equation}  \delta _{3} \leq  F(n) G(\infty ) ,  \Label{dFG}  \end{equation}
where
$$  G(\infty ) = \left[  2 - \left( \frac{4}{3} \right) ^{\beta _1} - \left( \frac{2}{3} \right) ^{\beta _1} \right]  .  $$
For later use, we define
$$  C_k = \left[  2 - \left( \frac{k+1}{k} \right) ^{\beta _1} - \left( \frac{k-1}{k} \right) ^{\beta _1} \right]  .  $$
Thus
\begin{eqnarray*}
   (n-3) \delta _{3} &\leq &  C_3 (n-3) \left( \frac{6n-9}{(n^2-1)(n-3)} \right) ^{\beta _1}  \\
   &=&  C_3  \left( \frac{(6n-9)(n-3)^{1/\beta _1-1}}{n^2-1} \right) ^{\beta _1}  \\
   &\leq &  \left . C_3  \left( \frac{(6n-9)(n-3)^{1/\beta _1-1}}{n^2-1} \right) ^{\beta _1}  \right |_{n=15} \\[1.2ex]
   &<& 0.0412 .   
\end{eqnarray*}
It is easy to verify that the expression in the second line is decreasing in $ n $
for $ n\geq 15 $, justifying the inequality given by third line. 

In the final step of the proof of Theorem~1, we need the sharper bound (\ref{dn6})
on $ (n-3)\delta _3 $. In the above proof, there is room for improvement
because we have used a fairly crude bound $ C_3 $ for the last
factor in (\ref{dn4}).

Using (\ref{dFG}) we see that (\ref{dn6}) is true for $ n\geq 99 $. It remains to
verify (\ref{dn6}) for the finite set $ n=45,\,47,\,\cdots ,\,97 $.
That can be achieved by brute force, simply
by computing each of these values numerically.
\end{PROOF}

\begin{LEMMA}
For fixed $ n\geq 45 $,
$$  a_{n-10}-a_{n-9}   < 0.1636.  $$
\end{LEMMA}  \begin{PROOF}
$$  a_{n-10}-a_{n-9} = \left( \frac{20n-100}{(n^2-1)(n-10)} \right) ^\beta  - \left( \frac{18n-91}{(n^2-1)(n-9)} \right) ^\beta    $$
has the form of $ \xi (\beta ) $ in (\ref{xi}), with $ \lambda =1 $, $ \mu =0 $, and satisfies Lemma~\ref{AlnA2} (ii).
This means that
we only need to find an upper bound when $ \beta =\beta _1 $. Thus
$$  (n-11)(a_{n-10}-a_{n-9}) \leq  \frac{(n-11)}{(n^2-1)^{\beta _1}}  \left[ \left( \frac{20n-100}{n-10} \right) ^{\beta _1} -\left( \frac{18n-81}{n-9} \right) ^{\beta _1} \right]   .  $$
The first factor is $ \dc $ in $ n\geq 45 $. The second factor is a difference
of two convex functions. The DIF technique proves that it is $ \dc $. Therefore,
\begin{eqnarray*}
   (n-11)(a_{n-10}-a_{n-9}) &\leq & \left . \frac{(n-11)}{(n^2-1)^{\beta _1}}  \left[ \left( \frac{20n-100}{n-10} \right) ^{\beta _1} -\left( \frac{18n-81}{n-9} \right) ^{\beta _1} \right]  \right |_{n=45} \\
   &<& 0.1636.
\end{eqnarray*}
\end{PROOF}

\end{document}